\documentclass[12pt,draft,leqno]{article}
\usepackage{amssymb, eucal, latexsym}

\newtheorem{theorem}{Theorem}[section]
\newtheorem{lm}[theorem]{Lemma}
\newtheorem{exa}[theorem]{Example}

\newtheorem{cor}[theorem]{Corollary}
\newtheorem{pro}[theorem]{Proposition}
\newtheorem{defi}[theorem]{Definition}

\newtheorem{nota}[theorem]{Notation}
\newtheorem{notas}[theorem]{Notations}

\newtheorem{nist}[theorem]{}

\def\p{\varphi}
\def\a{\alpha}
\def\b{\beta}
\def\d{\delta}

\def\g{\gamma}
\def\GA{\Gamma}

\def\l{\lambda}
\def\LAM{\Lambda}

\def\TE{\Theta}

\def\lra{\longrightarrow}

\def\sbe{\subseteq}
\def\spe{\supseteq}
\def\stm{\setminus}
\def\ems{\emptyset}
\def\nes{\neq\emptyset}

\def\ex{\exists}
\def\fa{\forall}

\def\we{\wedge}

\def\bv{\bigvee}

\def\ap{^\prime}
\def\inv{^{-1}}
\def\st{\ |\ }

\def\nin{\not\in}

\def\card #1{\vert #1 \vert}

\def\AA{{\cal A}}
\def\BB{{\cal B}}
\def\CC{{\cal C}}

\def\KK{{\cal K}}
\def\LL{{\cal L}}

\def\PP{{\cal P}}

\def\ZZ{{\cal Z}}

\def\ZLBA{{\bf ZLBA}}

\def\ZLC{{\bf ZLC}}

\def\ZLBA{{\bf ZLBA}}

\def\2{\mbox{{\bf 2}}}
\def\3{\mbox{{\bf 3}}}

\def\int{\mbox{{\rm int}}}

\def\cl{\mbox{{\rm cl}}}

\def\doc{\hspace{-1cm}{\em Proof.}~~}
\def\sq{\hspace*{\fill} \hbox{\vrule\vbox{\hrule\phantom{o}\hrule}\vrule}}
\def\sqs{\sq \vspace{2mm}}

\title{{\LARGE\bf
Open and other kinds of extensions over zero-dimensional local compactifications}\\
\vspace{0.35cm}
{\large\bf Georgi Dimov}\thanks{This paper was supported by the
project no. 005/2009 $``$General and Categorical Topology" of the
Sofia
University $``$St. Kl. Ohridski".}\\
\vspace{0.25cm}
 {\footnotesize Dept. of Math. and
Informatics, Sofia University,  Blvd. J. Bourchier 5, 1164 Sofia,
Bulgaria}}

\author{}

\date{}

\begin{document}
\maketitle
\begin{abstract}
{\footnotesize
\noindent  Generalizing a theorem of Ph. Dwinger \cite{Dw}, we
describe the partially ordered set of all (up to equivalence)
zero-dimensional locally compact Hausdorff extensions  of a
zero-dimensional Hausdorff space. Using this description, we find
the necessary and sufficient conditions which has to satisfy a map
between two zero-dimensional Hausdorff spaces in order to have
some kind of extension over arbitrary given in advance Hausdorff
zero-dimensional local compactifications of these spaces; we
regard the following kinds of extensions: continuous, open,
quasi-open, skeletal, perfect, injective, surjective. In this way
we generalize some classical results of B. Banaschewski \cite{Ba}
about the maximal zero-dimensional Hausdorff compactification.
Extending a recent theorem of G. Bezhanishvili \cite{B}, we
describe the local proximities corresponding to the
zero-dimensional Hausdorff local compactifications.}
\end{abstract}

{\footnotesize {\em  MSC:} primary 54C20, 54D35; secondary 54C10,
54D45, 54E05.

}

\footnotetext[1]{{\footnotesize {\em E-mail address:}
gdimov@fmi.uni-sofia.bg}}

\baselineskip = \normalbaselineskip

\section*{Introduction}

In \cite{Ba}, B. Banaschewski proved that every zero-dimensional
Hausdorff space $X$ has a zero-dimensional Hausdorff
compactification $\b_0X$ with the following remarkable property:
every continuous map $f:X\lra Y$, where $Y$  is a zero-dimensional
Hausdorff compact space, can be extended to a continuous map
$\b_0f:\b_0X\lra Y$; in particular, $\b_0X$ is the maximal
zero-dimensional Hausdorff compactification of $X$. As far as I
know, there are no descriptions of the maps $f$ for which the
extension $\b_0f$ is  open or quasi-open. In this paper we solve
the following more general problem: let $f:X\lra Y$ be a map
between two zero-dimensional Hausdorff spaces and $(lX,l_X)$,
$(lY,l_Y)$ be Hausdorff zero-dimensional locally compact
extensions of $X$ and $Y$, respectively; find the necessary and
sufficient conditions which has to satisfy the map $f$ in order to
have an $``$extension" $g:lX\lra lY$ (i.e. $g\circ l_X=l_Y\circ
f$) which is a map with some special properties (we regard the
following properties: continuous, open, perfect, quasi-open,
skeletal, injective, surjective). In \cite{LE2}, S. Leader solved
such a problem for continuous extensions over Hausdorff local
compactifications (= locally compact extensions)
 using the language of {\em local proximities} (the later,
as he showed, are in a bijective correspondence (preserving the
order) with the Hausdorff local compactifications regarded up to
equivalence). Hence, if one can describe the local proximities
which correspond to zero-dimensional Hausdorff local
compactifications then the above problem will be solved for
continuous extensions. Recently, G. Bezhanishvili \cite{B},
solving an old problem of L. Esakia, described the {\em
Efremovi\v{c} proximities} which correspond (in the sense of the
famous Smirnov Compactification Theorem \cite{Sm2}) to the
zero-dimensional Hausdorff compactifications
 (and called them {\em zero-dimensional Efremovi\v{c} proximities}).
 We extend here his result to the Leader's local proximities,
i.e. we describe the local proximities which correspond to the
Hausdorff zero-dimensional local compactifications and call them
{\em zero-dimensional local proximities} (see Theorem
\ref{zdlpth}). We do not use, however, these zero-dimensional
local proximities for solving our problem. We introduce a simpler
notion (namely, the {\em admissibe ZLB-algebra}) for doing this.
Ph. Dwinger \cite{Dw} proved, using Stone Duality Theorem
\cite{ST}, that the ordered set of all, up to equivalence,
zero-dimensional Hausdorff compactifications of a zero-dimensional
Hausdorff space is isomorphic to the ordered by inclusion set of
all {\em Boolean bases} of $X$ (i.e. of those bases of $X$ which
are Boolean subalgebras of the Boolean algebra $CO(X)$ of all
clopen (= closed and open) subsets of $X$). This description is
much simpler than that by Efremovi\v{c} proximities.  It was
rediscovered by K. D. Magill Jr. and J. A. Glasenapp \cite{MG} and
applied very successfully to the study of the poset of all, up to
equivalence, zero-dimensional Hausdorff compactifications of a
zero-dimensional Hausdorff space. We extend the cited above
Dwinger Theorem \cite{Dw} to the  zero-dimensional Hausdorff {\em
local compactifications} (see Theorem \ref{dwingerlc} below) with
the help of our generalization of the Stone Duality Theorem proved
in \cite{Di4} and the notion of $``$admissible ZLB-algebra" which
we introduce here. We obtain the solution of the problem
formulated above in the language of the admissible ZLB-algebras
(see Theorem \ref{zdextcmain}). As a corollary, we characterize
the maps $f:X\lra Y$ between two Hausdorff zero-dimensional spaces
$X$ and $Y$ for which the extension $\b_0f:\b_oX\lra\b_0Y$ is open
or quasi-open (see Corollary \ref{zdextcmaincb}).   Of course, one
can pass  from  admissible ZLB-algebras  to zero-dimensional local
proximities and conversely (see Theorem \ref{ailp} below; it
generalizes
  an analogous result about the connection between Boolean bases and
  zero-dimensional Efremovi\v{c} proximities
 obtained  in \cite{B}).

We now fix the notations.

 If $\CC$ denotes a category, we write
$X\in \card\CC$ if $X$ is
 an object of $\CC$, and $f\in \CC(X,Y)$ if $f$ is a morphism of
 $\CC$ with domain $X$ and codomain $Y$. By $Id_{\CC}$ we denote the identity functor on the category $\CC$.

All lattices are with top (= unit) and bottom (= zero) elements,
denoted respectively by 1 and 0. We do not require the elements
$0$ and $1$ to be distinct. Since we follow Johnstone's
terminology from \cite{J}, we will use the term {\em
pseudolattice} for a poset having all finite non-empty meets and
joins; the pseudolattices with a bottom will be called {\em
$\{0\}$-pseudolattices}. If $B$ is a Boolean algebra then we
denote by $Ult(B)$ the set of all ultrafilters in $B$.

If $X$ is a set then we denote the power set of $X$ by $P(X)$; the
identity function on $X$ is denoted by $id_X$.

 If
$(X,\tau)$ is a topological space and $M$ is a subset of $X$, we
denote by $\cl_{(X,\tau)}(M)$ (or simply by $\cl(M)$ or
$\cl_X(M)$) the closure of $M$ in $(X,\tau)$ and by
$\int_{(X,\tau)}(M)$ (or briefly by $\int(M)$ or $\int_X(M)$) the
interior of $M$ in $(X,\tau)$.

The  closed maps and the open maps between topological spaces are
assumed to be continuous but are not assumed to be onto. Recall
that a map is {\em perfect}\/ if it is closed and compact (i.e.
point inverses are compact sets).

For all notions and notations not defined here see \cite{Dw, E2, J, NW}.

\section{Preliminaries}

We will need some of our results from \cite{Di4} concerning the
extension of the Stone Duality Theorem to the category $\ZLC$ of
all locally compact zero-dimensio\-nal Hausdorff spaces and all
continuous maps between them.

Recall that if $(A,\le)$ is a poset and $B\sbe A$ then $B$ is said
to be a {\em dense subset of} $A$ if for any $a\in A\stm\{0\}$
there exists $b\in B\stm\{0\}$ such that $b\le a$.

\begin{defi}\label{deflba}{\rm (\cite{Di4})}
\rm A pair $(A,I)$, where $A$ is a Boolean algebra and $I$ is an
ideal of $A$ (possibly non proper) which is dense in $A$, is called a {\em local Boolean algebra} (abbreviated
as LBA).  Two LBAs $(A,I)$ and $(B,J)$ are said to be {\em
LBA-isomorphic} (or, simply, {\em isomorphic}) if there exists a
Boolean isomorphism $\p:A\lra B$ such that $\p(I)=J$.
\end{defi}

Let $A$ be a distributive $\{0\}$-pseudolattice and $Idl(A)$ be
the frame of all ideals of $A$. If $J\in Idl(A)$ then we will
write $\neg_A J$ (or simply $\neg J$) for the pseudocomplement of
$J$ in $Idl(A)$ (i.e. $\neg J=\bv\{I\in Idl(A)\st I\we
J=\{0\}\}$).  Recall that an ideal $J$ of $A$ is called {\em
simple} (Stone \cite{ST}) if $J\vee\neg J= A$ (i.e. $J$ has a
complement in $Idl(A)$). As it is proved in \cite{ST}, the set
$Si(A)$ of all simple ideals of $A$ is a Boolean algebra with
respect to the lattice operations in $Idl(A)$.

\begin{defi}\label{defzlba}{\rm (\cite{Di4})}
\rm An LBA $(B, I)$  is called a {\em ZLB-algebra} (briefly, {\em
ZLBA}) if, for every $J\in Si(I)$, the join $\bv_B J$($=\bv_B
\{a\st a\in J\}$) exists.

Let $\ZLBA$ be the category whose objects are all ZLBAs and whose
morphisms are all functions $\p:(B, I)\lra(B_1, I_1)$ between the
objects of $\ZLBA$ such that $\p:B\lra B_1$ is a Boolean
homomorphism satisfying the following condition:

\smallskip

\noindent(ZLBA) For every $b\in I_1$ there exists $a\in I$ such
that $b\le \p(a)$;

\smallskip

\noindent let the composition between the morphisms of $\ZLBA$ be
the usual composition between functions, and the
$\ZLBA$-identities be the identity functions.
\end{defi}

\begin{exa}\label{zlbaexa}{\rm (\cite{Di4})}
\rm Let $B$ be a Boolean algebra. Then the pair $(B,B)$ is a ZLBA.
\end{exa}

\begin{notas}\label{kxckx}
\rm Let $X$ be a topological space. We will denote by $CO(X)$ the
set of all clopen (= closed and open) subsets of $X$,
and by $CK(X)$ the set of all clopen compact subsets of $X$. For
every $x\in X$, we set
$u_x^{CO(X)}=\{F\in CO(X)\st x\in F\}.$
When there is no ambiguity, we will write $``u_x^C$" instead of
$``u_x^{CO(X)}$".
\end{notas}

The next assertion follows from the results obtained in
\cite{R2,Di4}.

\begin{pro}\label{psiult}
Let $(A,I)$ be a ZLBA. Set $X=\{u\in Ult(A)\st u\cap I\nes\}$.
Set, for every $a\in A$, $\l_A^C(a)=\{u\in X\st a\in u\}$. Let
$\tau$ be the topology on $X$ having as an open base the family
$\{\l_A^C(a)\st a\in I\}$. Then $(X,\tau)$ is a zero-dimensional
locally compact Hausdorff space, $\l_A^C(A)= CO(X)$,
$\l_A^C(I)=CK(X)$ and  $\l_A^C:A\lra CO(X)$ is a Boolean
isomorphism; hence, $\l_A^C:(A,I)\lra (CO(X),CK(X))$ is a
$\ZLBA$-isomorphism. We set $\TE^a(A,I)=(X,\tau)$.
\end{pro}

\begin{theorem}\label{genstonec}{\rm (\cite{Di4})}
The category\/ $\ZLC$  is dually equivalent to the category\/
$\ZLBA$.
In more details, let $\TE^a:\ZLBA\lra\ZLC$ and
$\TE^t:\ZLC\lra\ZLBA$ be two contravariant functors defined as
follows: for every $X\in\card{\ZLC}$, we set $\TE^t(X)=(CO(X),  CK(X))$,
and for every $f\in\ZLC(X,Y)$, $\TE^t(f):\TE^t(Y)\lra\TE^t(X)$ is
defined by the formula $\TE^t(f)(G)=f\inv(G)$, where $G\in CO(Y)$;
  for the definition of $\TE^a(B, I)$, where  $(B, I)$ is a ZLBA, see
  \ref{psiult};
 for every $\p\in\ZLBA((B, I),(B_1, J))$,
 $\TE^a(\p):\TE^a(B_1, J)\lra\TE^a(B, I)$ is given by the formula
$\TE^a(\p)(u\ap)=\p\inv(u\ap), \   \fa u\ap\in\TE^a(B_1,J)$; then
$t^{C}:Id_{\ZLC}\lra \TE^a\circ \TE^t$, where $t^{C}(X)=t_X^C, \
\fa X\in\card\ZLC$ and $t_X^C(x)=u_x^C$, for every $x\in X$, is a
natural isomorphism (hence, in particular,
$t_X^C:X\lra\TE^a(\TE^t(X))$ is a homeomorphism for every
$X\in\card{\ZLC}$); also, $\l^C: Id_{\ZLBA}\lra \TE^t\circ \TE^a$,
where $ \l^C(B, I)=\l_B^C, \ \fa (B, I)\in\card\ZLBA$, is a
natural isomorphism.
\end{theorem}

Finally, we will recall some definitions and facts from the theory
of extensions of topological spaces, as well as the fundamental
Leader's Local Compactification Theorem \cite{LE2}.

 Let $X$ be a Tychonoff space. We will denote by $\LL(X)$ the
set of all, up to equivalence, locally compact Hausdorff
extensions of $X$ (recall that two (locally compact Hausdorff)
extensions $(Y_1,f_1)$ and $(Y_2,f_2)$ of $X$ are said to be {\em
equivalent}\/ iff there  exists a homeomorphism $h:Y_1\lra Y_2$
such that $h\circ f_1=f_2$).  Let $[(Y_i,f_i)]\in\LL(X)$, where
$i=1,2$. We set $[(Y_1,f_1)]\le [(Y_2,f_2)]$  if there exists a
continuous  mapping $h:Y_2\lra Y_1$ such that $f_1=h\circ f_2$.
Then $(\LL(X),\le)$ is a poset (=partially ordered set).

 Let $X$ be a Tychonoff space. We will denote by $\KK(X)$ the
set of all, up to equivalence,  Hausdorff compactifications of
$X$.

\begin{nist}\label{nilea}
\rm Recall that if $X$ is a set and $P(X)$ is the power set of $X$
ordered by the inclusion, then a triple $(X,\d,\BB)$ is called
a {\em local proximity space} (see \cite{LE2}) if $\BB$ is an ideal (possibly non proper) of $P(X)$ and
$\d$ is
a symmetric binary relation
on $P(X)$ satisfying the following  conditions:

\smallskip

\noindent(P1) $\ems(-\d) A$ for every $A\sbe X$ ($``-\d$" means $``$not $\d$");

\smallskip

\noindent(P2) $A\d A$ for every $A\not =\ems$;

\smallskip

\noindent(P3) $A\d(B\cup C)$ iff $A\d B$ or $A\d C$;

\smallskip

\noindent(BC1) If $A\in \BB$, $C\sbe X$ and $A\ll C$ (where, for $D,E\sbe X$, $D\ll E$
iff $D(-\d) (X\stm E)$)
then there exists a $B\in\BB$
such that $A\ll B\ll C$;

\smallskip

\noindent(BC2) If $A\d C$, then there is a $B\in\BB$ such that $B\sbe C$ and $A\d B$.

\smallskip

\noindent A local proximity space $(X,\d,\BB)$ is said to be {\em
separated} if $\d$ is the identity relation on singletons. Recall
that every separated local proximity space $(X,\d,\BB)$ induces a
Tychonoff topology $\tau_{(X,\d,\BB)}$ in $X$ by defining
$\cl(M)=\{x\in X\st x\d M\}$ for every $M\sbe X$ (\cite{LE2}). If
$(X,\tau)$ is a topological space then we say that $(X,\d,\BB)$ is
a {\em local proximity space on} $(X,\tau)$ if
$\tau_{(X,\d,\BB)}=\tau$.

The  set of all separated local proximity spaces on a Tychonoff
space $(X,\tau)$ will be denoted by $\LL\PP(X,\tau)$. An order in
$\LL\PP(X,\tau)$ is defined by $(X,\b_1,\BB_1)\preceq
(X,\b_2,\BB_2)$ if $\b_2\sbe\b_1$ and $\BB_2\sbe\BB_1$ (see
\cite{LE2}).

A function $f:X_1\lra X_2$ between two local proximity spaces
$(X_1,\b_1,\BB_1)$
and
$(X_2,\b_2,\BB_2)$
is said to be an  {\em equicontinuous mapping\/} (see \cite{LE2}) if
the following two conditions are fulfilled:
\smallskip

\noindent(EQ1) $A\b_1 B$ implies $f(A)\b_2 f(B)$, for $A,B\sbe X$, and
\smallskip

\noindent(EQ2) $B\in\BB_1$ implies $f(B)\in\BB_2$.
\end{nist}

\begin{theorem}\label{Leader} {\rm (S. Leader \cite{LE2})}
Let $(X,\tau)$ be a Tychonoff space. Then there exists an
isomorphism $\LAM_X$ between the ordered sets $(\LL(X,\tau),\le)$
and $(\LL\PP(X,\tau),\preceq)$.  In more details, for every $(X,
\rho, \BB)\in\LL\PP(X,\tau)$ there exists a locally compact
Hausdorff extension $(Y,f)$  of X satisfying the following two
conditions:

\smallskip

\noindent(a) $A \rho B$ iff $\cl_Y(f(A))\cap \cl_Y(f(B))\nes$;

\smallskip

\noindent(b) $B\in\BB$  iff $\cl_Y(f(B))$ is compact.

\smallskip

\noindent Such a local compactification is unique up to
equivalence; we set $(Y,f)=L(X,\rho,\BB)$ and
$(\LAM_X)\inv(X,\rho,\BB)=[(Y,f)]$. The space $Y$ is compact iff
$X\in\BB$. Conversely, if $(Y,f)$ is a locally compact Hausdorff
extension of $X$ and $\rho$ and $\BB$ are defined by (a) and (b),
then $(X, \rho, \BB)$ is a separated local proximity space, and we
set $\LAM_X([(Y,f)])=(X,\rho,\BB)$.

Let $(X_i,\b_i,\BB_i)$, $i=1,2$, be two separated local proximity spaces and
$f:X_1\lra X_2$ be a function. Let $(Y_i,f_i)=L(X_i,\b_i,\BB_i)$, where $i=1,2$.
Then there exists a continuous map
$L(f):Y_1\lra Y_2$
such that $f_2\circ f= L(f)\circ f_1$ iff $f$ is an
equicontinuous map between
$(X_1,\b_1,\BB_1)$ and $(X_2,\b_2,\BB_2)$.
\end{theorem}

Recall that a subset $F$ of a topological space $(X,\tau)$ is
called {\em regular closed}\/ if $F=\cl(\int (F))$. Clearly, $F$
is regular closed iff it is the closure of an open set. For any
topological space $(X,\tau)$, the collection $RC(X,\tau)$ (we will
often write simply $RC(X)$) of all regular closed subsets of
$(X,\tau)$ becomes a complete Boolean algebra
$(RC(X,\tau),0,1,\we,\vee,{}^*)$ under the following operations: $
1 = X,  0 = \emptyset, F^* = \cl(X\stm F), F\vee G=F\cup G, F\we G
=\cl(\int(F\cap G)). $ The infinite operations are given by the
following formulas: $\bigvee\{F_\g\st
\g\in\GA\}=\cl(\bigcup\{F_\g\st \g\in\GA\})$ and
$\bigwedge\{F_\g\st \g\in\GA\}=\cl(\int(\bigcap\{F_\g\st
\g\in\GA\})).$ We denote by $CR(X,\tau)$ the family of all compact
regular closed subsets of $(X,\tau)$. We will often write  $CR(X)$
instead of $CR(X,\tau)$.

 We will  need a lemma from \cite{CNG}:

\begin{lm}\label{isombool}
Let $X$ be a dense subspace of a topological space $Y$. Then the
functions $r:RC(Y)\lra RC(X)$, $F\mapsto F\cap X$, and
$e:RC(X)\lra RC(Y)$, $G\mapsto \cl_Y(G)$, are Boolean isomorphisms
between Boolean algebras $RC(X)$ and $RC(Y)$, and $e\circ
r=id_{RC(Y)}$, $r\circ e=id_{RC(X)}$.
\end{lm}

\section{A Generalization of Dwinger Theorem}

\begin{defi}\label{admzlba}
\rm Let $X$ be a zero-dimensional Hausdorff space. Then:

\smallskip

\noindent(a) A ZLBA $(A,I)$ is called {\em admissible for} $X$ if
$A$ is a Boolean subalgebra of the Boolean algebra $CO(X)$ and $I$
is an open base of $X$.

\smallskip

\noindent(b) The set of all admissible  for $X$ ZLBAs is denoted by
$\ZZ\AA(X)$.

\smallskip

\noindent(c) If $(A_1,I_1),(A_2,I_2)\in \ZZ\AA(X)$  then we set
$(A_1,I_1)\preceq_0(A_2,I_2)$ if $A_1$ is a Boolean subalgebra of
$A_2$ and for every $V\in I_2$ there exists $U\in I_1$ such that
$V\sbe U$.
\end{defi}

\begin{nota}\label{lz}
\rm The set of all (up to equivalence) zero-dimensional locally
compact Hausdorff extensions of a zero-dimensional Hausdorff space
$X$ will be denoted by $\LL_0(X)$.
\end{nota}

\begin{theorem}\label{dwingerlc}
Let $X$ be a zero-dimensional Hausdorff space. Then the ordered
sets $(\LL_0(X),\le)$ and $(\ZZ\AA(X),\preceq_0)$ are isomorphic;
moreover, the zero-dimensional compact Hausdorff extensions of $X$
correspond to ZLBAs of the form $(A,A)$.
\end{theorem}

\doc
Let $(Y,f)$ be a locally compact  Hausdorff zero-dimensional
extensions of $X$. Set
\begin{equation}\label{01}
A_{(Y,f)}=f\inv(CO(Y)) \mbox{ and }
I_{(Y,f)}=f\inv(CK(Y)).
\end{equation}
 Note that $A_{(Y,f)}=\{F\in CO(X)\st
\cl_Y(f(F))$ is open in $Y\}$ and $I_{(Y,f)}=\{F\in A_{(Y,f)}\st
\cl_Y(f(F))$ is compact$\}$. We will show that
$(A_{(Y,f)},I_{(Y,f)})\in\ZZ\AA(X)$. Obviously, the map
$r_{(Y,f)}^0:(CO(Y),CK(Y))\lra(A_{(Y,f)},I_{(Y,f)}), \ G\mapsto
f\inv(G),$
is a Boolean isomorphism  such that
$r_{(Y,f)}^0(CK(Y))=I_{(Y,f)}$. Hence $(A_{(Y,f)},I_{(Y,f)})$ is a
ZLBA and $r_{(Y,f)}^0$ is an LBA-isomorphism.  It is easy to see
that $I_{(Y,f)}$ is a base of $X$ (because $Y$ is locally
compact). Hence $(A_{(Y,f)},I_{(Y,f)})\in\ZZ\AA(X)$.
 It is
clear that  if $(Y_1,f_1)$ is a locally compact Hausdorff
zero-dimensional extensions of $X$ equivalent to the extension
$(Y,f)$, then
$(A_{(Y,f)},I_{(Y,f)})=(A_{(Y_1,f_1)},I_{(Y_1,f_1)})$. Therefore,
a map
\begin{equation}\label{dw1}
\a_X^0:\LL_0(X)\lra\ZZ\AA(X), \
[(Y,f)]\mapsto(A_{(Y,f)},I_{(Y,f)}),
\end{equation}
 is well-defined.

 Let $(A,I)\in\ZZ\AA(X)$ and
$Y=\TE^a(A,I)$. Then $Y$ is a locally compact Hausdorff
zero-dimensional space. For every $x\in X$, set
\begin{equation}\label{02}
u_{x,A}=\{F\in
A\st x\in F\}.
\end{equation}
 Since $I$ is a base of $X$, we get that $u_{x,A}$
is an ultrafilter in $A$ and $u_{x,A}\cap I\nes$, i.e. $u_{x,A}\in
Y$. Define
\begin{equation}\label{0f}
f_{(A,I)}:X\lra Y, \   x\mapsto u_{x,A}.
\end{equation}
 Set, for
short, $f=f_{(A,I)}$. Obviously, $\cl_Y(f(X))=Y$.  It is easy to
see that $f$ is a homeomorphic embedding. Hence $(Y,f)$ is a
locally compact Hausdorff zero-dimensional extension of $X$. We
now set:
\begin{equation}\label{dw2}
\b_X^0:\ZZ\AA(X)\lra \LL_0(X), \  (A,I)\mapsto
[(\TE^a(A,I),f_{(A,I)})].
\end{equation}

We will show that $\a_X^0\circ\b_X^0=id_{\ZZ\AA(X)}$ and
$\b_X^0\circ\a_X^0=id_{\LL_0(X)}$.

Let $[(Y,f)]\in\LL_0(X)$. Set, for short, $A=A_{(Y,f)}$,
$I=I_{(Y,f)}$, $g=f_{(A,I)}$, $Z=\TE^a(A,I)$ and $\p=r_{(Y,f)}^0$.
Then $\b_X^0(\a_X^0([(Y,f)]))=\b_X(A,I)= [(Z,g)]$. We have to show
that $[(Y,f)]=[(Z,g)]$. Since $\p$ is an LBA-isomorphism, we get
that $h=\TE^a(\p):Z\lra\TE^a(\TE^t(Y))$ is a homeomorphism. Set
$Y\ap=\TE^a(\TE^t(Y))$.  By  Theorem \ref{genstonec}, the map
$t_Y^C:Y\lra Y\ap, \ y\mapsto u_y^{CO(Y)}$ is a homeomorphism. Let
$h\ap=(t_Y^C)\inv\circ h$. Then $h\ap:Z\lra Y$ is a homeomorphism.
We will prove that $h\ap\circ g=f$ and this will imply that
$[(Y,f)]=[(Z,g)]$. Let $x\in X$. Then
$h\ap(g(x))=h\ap(u_{x,A})=(t_Y^C)\inv(h(u_{x,A}))=(t_Y)\inv(\p\inv(u_{x,A}))$.
We have that $u_{x,A}=\{f\inv(F)\st F\in CO(Y), x\in
f\inv(F)\}=\{\p(F)\st F\in CO(Y), x\in f\inv(F)\}$. Thus
$\p\inv(u_{x,A})=\{F\in CO(Y)\st f(x)\in F\}=u_{f(x)}^{CO(Y)}$. Hence
$(t_Y)\inv(\p\inv(u_{x,A}))=f(x)$. So, $h\ap\circ g=f$. Therefore,
$\b_X^0\circ\a_X^0=id_{\LL_0(X)}$.

Let $(A,I)\in\ZZ\AA(X)$ and $Y=\TE^a(A,I)$.  Set $f=f_{(A,I)}$,
$B=A_{(Y,f)}$ and $J=I_{(Y,f)}$. Then
 $\a_X^0(\b_X^0(A,I))=(B,J)$.
   By  Theorem \ref{genstonec}, we have that
$\l_A^C:(A,I)\lra(CO(Y),CK(Y))$ is an LBA-isomorphism. Hence
$\l_A^C(A)=CO(Y)$ and $\l_A^C(I)=CK(Y)$. We will show that
$f\inv(\l_A^C(F))=F$, for every $F\in A$. Recall that
$\l_A^C(F)=\{u\in Y\st F\in u\}$.  Now we have that if $F\in A$
then $f\inv(\l_A^C(F))=\{x\in X\st f(x)\in\l_A^C(F)\}=\{x\in X\st
u_{x,A}\in\l_A^C(F)\}=\{x\in X\st F\in u_{x,A}\}=\{x\in X\st x\in
F\}=F$. Thus
\begin{equation}\label{05}
B=f\inv(CO(Y))=A \mbox{ and } J=f\inv(CK(Y))=I.
\end{equation}
 Therefore,
$\a_X^0\circ\b_X^0=id_{\ZZ\AA(X)}$.

We will now prove that $\a_X^0$ and $\b_X^0$ are monotone maps.

Let $[(Y_i,f_i)]\in\LL_0(X)$, where $i=1,2$, and
$[(Y_1,f_1)]\le[(Y_2,f_2)]$. Then there exists a continuous map
$g:Y_2\lra Y_1$ such that $g\circ f_2=f_1$. Set
$A_i=A_{(Y_i,f_i)}$ and $I_i=I_{(Y_i,f_i)}$, $i=1,2$. Then
$\a_X^0([(Y_i,f_i)])=(A_i,I_i)$, where $i=1,2$. We have to show
that $A_1\sbe A_2$ and  for every $V\in I_2$ there exists $U\in
I_1$ such that $V\sbe U$. Let $F\in A_1$. Then
$F\ap=\cl_{Y_1}(f_1(F))\in CO(Y_1)$ and, hence, $G\ap=g\inv(F\ap)\in
CO(Y_2)$. Thus $(f_2)\inv(G\ap)\in A_2$. Since
$(f_2)\inv(G\ap)=(f_2)\inv(g\inv(F\ap))=(f_2)\inv(g\inv(\cl_{Y_1}(f_1(F))))=(f_1)\inv(\cl_{Y_1}(f_1(F)))=F$,
we get that $F\in A_2$. Therefore, $A_1\sbe A_2$. Further, let
$V\in I_2$. Then $V\ap=\cl_{Y_2}(f_2(V))\in CK(Y_2)$. Thus $g(V\ap)$ is
a compact subset of $Y_1$. Hence there exists $U\in I_1$ such that
$g(V\ap)\sbe \cl_{Y_1}(f_1(U))$. Then $V\sbe
(f_2)\inv(g\inv(g(\cl_{Y_2}(f_2(V)))))=(f_1)\inv(g(V\ap))\sbe
(f_1)\inv(\cl_{Y_1}(f_1(U)))=U$. So,
$\a_X^0([(Y_1,f_1)])\preceq_0\a_X^0([(Y_2,f_2)])$. Hence, $\a_X^0$
is a monotone function.

Let now $(A_i,I_i)\in\ZZ\AA(X)$, where $i=1,2$, and
$(A_1,I_1)\preceq_0(A_2,I_2)$. Set, for short,
$Y_i=\TE^a(A_i,I_i)$ and $f_i=f_{(A_i,I_i)}$, $i=1,2$. Then
$\b_X^0(A_i,I_i)=[(Y_i,f_i)]$, $i=1,2$. We will show that
$[(Y_1,f_1)]\le[(Y_2,f_2)]$. We have that, for $i=1,2$, $f_i:X\lra
Y_i$ is defined by $f_i(x)=u_{x,A_i}$, for every $x\in X$. We also
have that $A_1\sbe A_2$ and  for every $V\in I_2$ there exists
$U\in I_1$ such that $V\sbe U$. Let us regard the function
$\p:(A_1,I_1)\lra(A_2,I_2), \ F\mapsto F.$
Obviously, $\p$ is a $\ZLBA$-morphism.  Then $g=\TE^a(\p):Y_2\lra
Y_1$ is a continuous map. We will prove that $g\circ f_2=f_1$,
i.e. that for every $x\in X$, $g(u_{x,A_2})=u_{x,A_1}$. So, let
$x\in X$. We have that $u_{x,A_2}=\{F\in A_2\st x\in F\}$ and
$g(u_{x,A_2})=\p\inv(u_{x,A_2})$. Clearly,
$\p\inv(u_{x,A_2})=\{F\in A_1\cap A_2\st x\in F\}$. Since $A_1\sbe
A_2$, we get that $\p\inv(u_{x,A_2})=\{F\in A_1\st x\in
F\}=u_{x,A_1}$. So, $g\circ f_2=f_1$. Thus
$[(Y_1,f_1)]\le[(Y_2,f_2)]$. Therefore, $\b_X^0$ is also a
monotone function. Since $\b_X^0=(\a_X^0)\inv$, we get that
$\a_X^0$ (as well as $\b_X^0$) is an isomorphism. \sqs

\begin{defi}\label{admba}
\rm Let $X$ be a zero-dimensional Hausdorff space.
 A Boolean algebra $A$ is called {\em admissible for}
$X$ (or, a {\em Boolean base of} $X$) if $A$ is a Boolean
subalgebra of the Boolean algebra $CO(X)$ and $A$ is an open base
of $X$.
 The set of all admissible Boolean algebras for $X$ is
denoted by $\BB\AA(X)$.
\end{defi}

\begin{nota}\label{cz}
\rm The set of all (up to equivalence) zero-dimensional compact
Hausdorff
  extensions of a zero-dimensional Hausdorff space $X$
will be denoted by $\KK_0(X)$.
\end{nota}

\begin{cor}\label{dwinger}{\rm (Ph. Dwinger \cite{Dw})}
Let $X$ be a zero-dimensional Hausdorff space. Then the ordered
sets $(\KK_0(X),\le)$ and $(\BB\AA(X),\sbe)$ are isomorphic.
\end{cor}

\doc Clearly, a  Boolean algebra $A$ is  admissible for
$X$ iff the ZLBA $(A,A)$ is  admissible for $X$. Also, if
$A_1,A_2$ are two   admissible for $X$ Boolean algebras then
$A_1\sbe A_2$ iff $(A_1,A_1)\preceq_0(A_2,A_2)$. Since the
admissible ZLBAs of the form $(A,A)$ and only they correspond to
the zero-dimensional compact Hausdorff extensions of $X$, it
becomes obvious that our assertion follows from Theorem
\ref{dwingerlc}. \sqs

\section{Zero-dimensional Local Proximities}

\begin{defi}\label{zdlpdef}
\rm A local proximity $(X,\d,\BB)$ is called {\em zero-dimensional}
if for every $A,B\in\BB$ with $A\ll B$ there exists $C\sbe X$
such that $A\sbe C\sbe B$ and $C\ll C$.

The  set of all separated zero-dimensional local proximity spaces on a Tychonoff
space $(X,\tau)$ will be denoted by $\LL\PP_0(X,\tau)$. The restriction of the
order relation $\preceq$ in $\LL\PP(X,\tau)$ (see \ref{nilea}) to the set
$\LL\PP_0(X,\tau)$ will be denoted again by $\preceq$.
\end{defi}

\begin{theorem}\label{zdlpth}
Let $(X,\tau)$ be a zero-dimensional Hausdorff space. Then the ordered
sets $(\LL_0(X),\le)$ and $(\LL\PP_0(X,\tau),\preceq)$ are isomorphic
(see \ref{zdlpdef} and \ref{dwingerlc} for the notations).
\end{theorem}

\doc Having in mind Leader's Theorem \ref{Leader}, we need only to show that if
$[(Y,f)]\in \LL(X)$ and $\LAM_X([(Y,f)])=(X,\d,\BB)$ then $Y$ is a
zero-dimensional space iff $(X,\d,\BB)\in\LL\PP_0(X)$.

So, let $Y$ be a zero-dimensional space. Then, by Theorem
\ref{Leader}, $\BB=\{B\sbe X\st \cl_Y(f(B))$ is compact$\}$, and
for every $A,B\sbe X$, $A\d B$ iff
$\cl_Y(f(A))\cap\cl_Y(f(B))\nes$. Let $A,B\in\BB$ and $A\ll B$.
Then $\cl_Y(f(A))\cap\cl_Y(f(X\stm B))=\ems$. Since $\cl_Y(f(A))$
is compact and $Y$ is zero-dimensional, there exists $U\in CO(Y)$
such that $\cl_Y(f(A))\sbe U\sbe Y\stm \cl_Y(f(X\stm B))$. Set
$V=f\inv(U)$. Then $A\sbe V\sbe\int_X(B)$, $\cl_Y(f(V))=U$ and
$\cl_Y(f(X\stm V))=Y\stm U$. Thus $V\ll V$ and $A\sbe V\sbe B$.
 Therefore, $(X,\d,\BB)\in\LL\PP_0(X)$.

Conversely, let $(X,\d,\BB)\in\LL\PP_0(X)$ and $(Y,f)=L(X,\d,\BB)$
(see \ref{Leader} for the notations). We will prove that $Y$ is a
zero-dimensional space. We have again, by Theorem \ref{Leader},
that the
 formulas written in the preceding paragraph for $\BB$ and $\d$ take place.
 Let $y\in Y$ and $U$ be an open neighborhood of $y$.
Since $Y$ is locally compact and Hausdorff, there exist $F_1, F_2\in CR(Y)$
such that $y\in F_1\sbe\int_Y(F_2)\sbe F_2\sbe U$. Let $A_i=f\inv(F_i)$, $i=1,2$.
Then $\cl_Y(f(A_i))=F_i$,  and hence $A_i\in\BB$, for $i=1,2$. Also, $A_1\ll A_2$.
Thus there exists $C\in\BB$ such that $A_1\sbe C\sbe A_2$ and $C\ll C$.
It is easy to see that $F_1\sbe\cl_Y(f(C))\sbe F_2$ and that $\cl_Y(f(C))\in CO(Y)$.
Therefore, $Y$ is a zero-dimensional space. \sqs

By Theorem \ref{Leader}, for every Tychonoff space $(X,\tau)$, the
local proximities of the form $(X,\d,P(X))$ on $(X,\tau)$ and only
they correspond to the Hausdorff compactifications of $(X,\tau)$.
The pairs $(X,\d)$ for which the triple $(X,\d,P(X))$ is a local
proximity are called {\em Efremovi\v{c} proximities}. Hence,
Leader's Theorem \ref{Leader} implies  the famous Smirnov
Compactification Theorem \cite{Sm2}. The notion of a
zero-dimensional proximity was introduced recently by G.
Bezhanishvili \cite{B}. Our notion of  a zero-dimensional local
proximity is a generalization of it. We will denote by $\PP_0(X)$
the set of all  zero-dimensional proximities on a zero-dimensional
Hausdorff  space  $X$. Now it becomes clear that our Theorem
\ref{zdlpth} implies immediately the following theorem of G.
Bezhanishvili \cite{B}:

\begin{cor}\label{zdlpcor}{\rm (G. Bezhanishvili \cite{B})}
Let $(X,\tau)$ be a zero-dimensional Hausdorff space.
Then there exists an isomorphism between the ordered
sets $(\KK_0(X),\le)$ and $(\PP_0(X,\tau),\preceq)$
(see \ref{zdlpdef} and \ref{dwingerlc} for the notations).
\end{cor}

The connection between the zero-dimensional local proximities on a
zero-dimensional Hausdorff space $X$ and
the admissible for $X$ ZLBAs is clarified in the next result:

\begin{theorem}\label{ailp}
Let $(X,\tau)$ be a zero-dimensional Hausdorff space. Then:

\smallskip

\noindent(a) Let $(A,I)\in\ZZ\AA(X,\tau)$. Set $\BB=\{M\sbe X\st \ex B\in I$ such that $M\sbe B\}$,
and for every $M,N\in\BB$, let $M\d N\iff (\fa F\in I)[(M\sbe F)\rightarrow(F\cap N\nes)]$; further, for every
$K,L\sbe X$, let $K\d L\iff[\ex M,N\in\BB$ such that $M\sbe K, N\sbe L$ and $M\d N]$. Then $(X,\d,\BB)\in\LL\PP_0(X,\tau)$.
Set $(X,\d,\BB)=L_X(A,I)$.

\smallskip

\noindent(b) Let $(X,\d,\BB)\in\LL\PP_0(X,\tau)$. Set $A=\{F\sbe X\st F\ll F\}$ and $I=A\cap\BB$. Then $(A,I)\in\ZZ\AA(X,\tau)$.
Set $(A,I)=l_X(X,\d,\BB)$.

\smallskip

\noindent(c)  $\b_X^0=(\LAM_X)\inv\circ L_X$ and, for every
$(X,\d,\BB)\in\LL\PP_0(X,\tau)$, $(\b_X^0\circ
l_X)(X,\d,\BB)=(\LAM_X)\inv(X,\d,\BB)$ (see \ref{Leader},
(\ref{dw2}), as well as (a) and (b) here for the notations);

\smallskip

\noindent(d) The correspondence
$L_X:(\ZZ\AA(X,\tau),\preceq_0)\lra(\LL\PP_0(X,\tau),\preceq)$ is
an isomorphism (between posets) and $L_X\inv=l_X$.
\end{theorem}

\doc It follows from Theorems \ref{dwingerlc}, \ref{zdlpth} and \ref{Leader}. \sqs

The above assertion is a generalization of the analogous result of G. Bezhanishvili \cite{B} concerning the
connection between the zero-dimensional  proximities on a zero-dimensional Hausdorff space $X$ and
the admissible for $X$ Boolean algebras.

\section{Extensions over Zero-dimensional Local Compactifications}

\begin{theorem}\label{zdextc}
Let $(X_i,\tau_i)$, where $i=1,2$, be  zero-dimensional Hausdorff
spaces, $(Y_i,f_i)$ be a zero-dimensional Hausdorff local
compactification of $(X_i,\tau_i)$, $(A_i,I_i)=\a_X^0(Y_i,f_i)$
 (see (\ref{dw1}) and (\ref{01})  for
$\a_{X_i}^0$)),  where $i=1,2$,
 and $f:X_1\lra X_2$ be a function. Then there exists a continuous function
 $g=L_0(f):Y_1\lra Y_2$ such that $g\circ f_1=f_2\circ f$ iff $f$ satisfies the
 following conditions:

\smallskip

\noindent{\rm (ZEQ1)} For every $G\in A_2$, $f\inv(G)\in A_1$ holds;

\smallskip

\noindent{\rm (ZEQ2)} For every $F\in I_1$ there exists $G\in I_2$ such that $f(F)\sbe G$.
\end{theorem}

\doc ($\Rightarrow$) Let there exists a continuous function $g:Y_1\lra Y_2$
such that $g\circ f_1=f_2\circ f$. By Lemma \ref{isombool} and (\ref{05}),
we have that the maps
\begin{equation}\label{0i}
r_i^c:CO(Y_i)\lra A_i, \ G\mapsto (f_i)\inv(G), \ e_i^c:A_i\lra CO(Y_i), \ F\mapsto\cl_{Y_i}(f_i(F)),
\end{equation}
where $i=1,2$, are Boolean isomorphisms; moreover, since  $r_i^c(CK(Y_i))=I_i$
and $e_i^c(I_i)=CK(Y_i)$, we get that
\begin{equation}\label{0il}
r_i^c:(CO(Y_i),CK(Y_i))\lra (A_i,I_i) \mbox{ and }
e_i^c:(A_i,I_i)\lra (CO(Y_i),CK(Y_i)),
\end{equation}
 where $i=1,2$, are LBA-isomorphisms.
Set
\begin{equation}\label{0ic}
\psi_g:CO(Y_2)\lra CO(Y_1), \ G\mapsto g\inv(G), \mbox{ and }
\psi_f=r_1^c\circ\psi_g\circ e_2^c.
\end{equation}
Then $\psi_f:A_2\lra A_1$. We will prove that
\begin{equation}\label{psif}
\psi_f(G)=f\inv(G), \mbox{ for every } G\in A_2.
\end{equation}
Indeed, let $G\in A_2$.
Then $\psi_f(G)=(r_1^c\circ\psi_g\circ e_2^c)(G)=(f_1)\inv(g\inv(\cl_{Y_2}(f_2(G))))=
\{x\in X_1\st (g\circ f_1)(x)\in\cl_{Y_2}(f_2(G))\}=\{x\in X_1\st f_2(f(x))\in\cl_{Y_2}(f_2(G))\}
=\{x\in X_1\st f(x)\in (f_2)\inv(\cl_{Y_2}(f_2(G)))\}=
\{x\in X_1\st f(x)\in G\}=f\inv(G)$. This shows that condition (ZEQ1) is fulfilled.
Since, by Theorem \ref{genstonec}, $\psi_g=\TE^t(g)$, we get that $\psi_g$ is
a $\ZLBA$-morphism. Thus $\psi_f$ is a $\ZLBA$-morphism. Therefore, for every $F\in I_1$
there exists $G\in I_2$ such that $f\inv(G)\spe F$. Hence, condition (ZEQ2)
is also checked.

\smallskip

\noindent($\Leftarrow$) Let $f$ be a function satisfying
conditions (ZEQ1) and (ZEQ2). Set $\psi_f:A_2\lra A_1$, $G\mapsto
f\inv(G)$. Then $\psi_f:(A_2,I_2)\lra (A_1,I_1)$ is a
$\ZLBA$-morphism. Put $g=\TE^a(\psi_f)$. Then
$g:\TE^a(A_1,I_1)\lra\TE^a(A_2,I_2)$, i.e. $g:Y_1\lra Y_2$ and $g$
is a continuous function (see Theorem \ref{genstonec} and
(\ref{dw2})). We will show that $g\circ f_1=f_2\circ f$. Let $x\in
X_1$. Then, by (\ref{0f}) and Theorem \ref{genstonec}, $g(f_1(x))=
g(u_{x,A_1})=(\psi_f)\inv(u_{x,A_1})=\{G\in A_2\st\psi_f(G)\in
u_{x,A_1}\}=\{G\in A_2\st x\in f\inv(G)\}=\{G\in A_2\st f(x)\in
G\}=u_{f(x),A_2}=f_2(f(x))$. Thus, $g\circ f_1=f_2\circ f$. \sqs

It is natural to write $f:(X_1,A_1,I_1)\lra (X_2,A_2,I_2)$ when we
have a situation like that which is described in Theorem
\ref{zdextc}. Then, in analogy with the Leader's equicontinuous
functions (see Leader's Theorem \ref{Leader}), the functions
$f:(X_1,A_1,I_1)\lra (X_2,A_2,I_2)$ which satisfy conditions
(ZEQ1) and (ZEQ2) will be called {\em 0-equicontinuous functions}.
Since $I_2$ is a base of $X_2$,
 we obtain that every 0-equcontinuous function is a continuous function.

\begin{cor}\label{zdextcc}
Let $(X_i,\tau_i)$, $i=1,2$, be two zero-dimensional Hausdorff spaces,
$A_i\in\BB\AA(X_i)$, $(Y_i,f_i)=\b_{X_i}^0(A_i,A_i)$ (see (\ref{dw2}) for $\b_{X_i}^0$),
where $i=1,2$,
 and $f:X_1\lra X_2$ be a function. Then there exists a continuous function
 $g=L_0(f):Y_1\lra Y_2$ such that $g\circ f_1=f_2\circ f$ iff $f$ satisfies  condition
{\rm (ZEQ1)}.
\end{cor}

\doc It follows from Theorem \ref{zdextc} because for ZLBAs of the form $(A_i,A_i)$,
where $i=1,2$, condition (ZEQ2) is always fulfilled. \sqs

Clearly, Theorem \ref{dwinger} implies (and this is noted in
\cite{Dw}) that every zero-dimensi\-o\-nal Hausdorff space $X$ has
a greatest zero-dimensional Hausdorff compactification which
corresponds to the admissible for $X$ Boolean algebra $CO(X)$.
This compactification was discovered by B. Banaschewski \cite{Ba};
it is denoted by $(\b_0X,\b_0)$ and it is called {\em Banaschewski
compactification} of $X$. One obtains immediately its main
property using our Corollary \ref{zdextcc}:

\begin{cor}\label{zdextcb}{\rm (B. Banaschewski \cite{Ba})}
Let $(X_i,\tau_i)$, $i=1,2$, be two zero-dimensional Hausdorff spaces and
$(cX_2,c)$ be a zero-dimensional Hausdorff compactification of $X_2$. Then for every continuous
function $f:X_1\lra X_2$ there exists a continuous function
$g:\b_0X_1\lra cX_2$ such that $g\circ\b_0=c\circ f$.
\end{cor}

\doc Since $\b_0X_1$ corresponds to the admissible for $X_1$ Boolean algebra $CO(X_1)$,
condition (ZEQ1) is clearly fulfilled when $f$ is a continuous function. Now apply
Corollary \ref{zdextcc}. \sqs

If in the above Corollary \ref{zdextcb} $cX_2=\b_0X_2$, then the map $g$ will be denoted by $\b_0f$.

Recall
that a function $f:X\lra Y$ is called {\em skeletal}\/ (\cite{MR})
if
\begin{equation}\label{ske}
\int(f\inv(\cl (V)))\sbe\cl(f\inv(V))
\end{equation}
for every open subset  $V$  of $Y$. Recall also the following result:

\begin{lm}\label{skel}{\rm (\cite{D1})}
A function $f:X\lra Y$ is  skeletal iff\/ $\int(\cl(f(U)))\nes$,
for every  non-empty  open subset $U$ of $X$.
\end{lm}

\begin{lm}\label{skelnew}
A  continuous map $f:X\lra Y$, where  $X$ and $Y$ are topological spaces, is skeletal iff
for every open subset $V$ of $Y$ such that $\cl_Y(V)$ is open,
$\cl_X(f\inv(V))= f\inv(\cl_Y(V))$ holds.
\end{lm}

\doc ($\Rightarrow$) Let $f$ be a skeletal continuous map and $V$ be
an open subset of $Y$ such that $\cl_Y(V)$ is open. Let $x\in
f\inv(\cl_Y(V))$. Then $f(x)\in \cl_Y(V)$. Since $f$ is
continuous, there exists an open neighborhood $U$ of $x$ in $X$
such that $f(U)\sbe \cl_Y(V)$. Suppose that
$x\nin\cl_X(f\inv(V))$. Then there exists an open neighborhood $W$
of $x$ in $X$ such that $W\sbe U$ and $W\cap f\inv(V)=\ems$. We
obtain that
 $\cl_Y(f(W))\cap V=\ems$ and $\cl_Y(f(W))\sbe \cl_Y(f(U))\sbe
\cl_Y(V)$. Since, by Lemma \ref{skel}, $\int_Y(\cl_Y(f(W)))\nes$, we get a contradiction.
Thus $f\inv(\cl_Y(V))\sbe\cl_X(f\inv(V))$.
The converse inclusion follows from the continuity of $f$.
Hence $f\inv(\cl_Y(V))=\cl_X(f\inv(V))$.

\smallskip

\noindent($\Leftarrow$) Suppose that there exists an open  subset
$U$ of $X$ such that $\int_Y(\cl_Y(f(U)))=\ems$ and $U\nes$. Then,
clearly, $V=Y\stm\cl_Y(f(U))$ is an open dense subset of $Y$.
Hence $\cl_Y(V)$ is open in $Y$. Thus $\cl_X(f\inv(V))=
f\inv(\cl_Y(V))=f\inv(Y)=X$ holds. Therefore
$X=\cl_X(f\inv(V))=\cl_X(f\inv(Y\stm \cl_Y(f(U))))=\cl_X(X\stm
f\inv(\cl_Y(f(U))))$. Since $U\sbe f\inv(\cl_Y(f(U)))$, we get
that $X\stm U\spe\cl_X(X\stm f\inv(\cl_Y(f(U))))=X$, a
contradiction. Hence, $f$ is a skeletal map. \sqs

Note that the proof of Lemma \ref{skelnew} shows that the following assertion is also true:

\begin{lm}\label{skelnewcor}
A  continuous map $f:X\lra Y$, where  $X$ and $Y$ are topological spaces, is skeletal iff
for every open dense subset $V$ of $Y$, $\cl_X(f\inv(V))= X$ holds.
\end{lm}

\begin{lm}\label{skelnewnew}
Let $(X_i,\tau_i)$, $i=1,2$, be two topological spaces,
$(Y_i,f_i)$ be some extensions of $(X_i,\tau_i)$, $i=1,2$,
$f:X_1\lra X_2$ and $g:Y_1\lra Y_2$ be two continuous functions
such that $g\circ f_1=f_2\circ f$. Then $g$ is skeletal iff $f$ is
skeletal.
\end{lm}

\doc ($\Rightarrow$) Let $g$ be skeletal and $V$ be an open dense
subset of $X_2$. Set $U=Ex_{Y_2}(V)$, i.e. $U=Y_2\stm
\cl_{Y_2}(f_2(X_2\stm V))$. Then $U$ is an open dense subset of
$Y_2$ and $f_2\inv(U)=V$. Hence, by Lemma \ref{skelnewcor},
$g\inv(U)$ is a dense open subset of $Y_1$. We will prove that
$f_1\inv(g\inv(U))\sbe f\inv(V)$. Indeed, let $x\in
f_1\inv(g\inv(U))$. Then $g(f_1(x))\in U$, i.e. $f_2(f(x))\in U$.
Thus $f(x)\in f_2\inv(U)=V$. So, $f_1\inv(g\inv(U))\sbe f\inv(V)$.
This shows that $f\inv(V)$ is dense in $X_1$. Therefore, by Lemma
\ref{skelnewcor}, $f$ is a skeletal map.

\smallskip

\noindent($\Leftarrow$) Let $f$ be a skeletal map and $U$ be a
dense open subset of $Y_2$. Set $V=f_2\inv(U)$. Then $V$ is an
open dense subset of $X_2$. Thus, by Lemma \ref{skelnewcor},
$f\inv(V)$ is a dense subset of $X_1$. We will prove that
$f\inv(V)\sbe f_1\inv(g\inv(U))$. Indeed, let $x\in f\inv(V)$.
Then $f(x)\in V=f_2\inv(U)$. Thus $f_2(f(x))\in U$, i.e.
$g(f_1(x))\in U$. So, $f\inv(V)\sbe f_1\inv(g\inv(U))$. This
implies that $g\inv(U)$ is dense in $Y_1$. Now, Lemma
\ref{skelnewcor} shows that $g$ is a skeletal map. \sqs

We are now ready to prove the following result:

\begin{theorem}\label{zdextcmain}
Let $(X_i,\tau_i)$, where $i=1,2$, be  zero-dimensional Hausdorff
spaces. Let, for $i=1,2$,  $(Y_i,f_i)$ be a zero-dimensional
Hausdorff local compactification of $(X_i,\tau_i)$,
$(A_i,I_i)=\a_X^0(Y_i,f_i)$
 (see (\ref{dw1}) and (\ref{01})  for
$\a_{X_i}^0$),
  $f:(X_1,A_1,I_1)\lra (X_2,A_2,I_2)$
   be a 0-equicontinuous function and  $g=L_0(f):Y_1\lra Y_2$ be the continuous
   function such that $g\circ f_1=f_2\circ f$ (its existence is guaranteed by
  Theorem \ref{zdextc}). Then:

  \smallskip

\noindent(a) $g$ is skeletal iff $f$ is skeletal;

  \smallskip

\noindent(b) $g$ is an open map iff $f$ satisfies the following condition:

\smallskip

\noindent{\rm(ZO)} For every $F\in I_1$,  $\cl_{X_2}(f(F))\in I_2$ holds;

  \smallskip

\noindent(c) $g$ is a perfect map iff $f$ satisfies the following condition:

\smallskip

\noindent{\rm(ZP)} For every $G\in I_2$,  $f\inv(G)\in I_1$ holds (i.e., briefly, $f\inv(I_2)\sbe I_1$);

\smallskip

\noindent(d) $\cl_{Y_2}(g(Y_1))=Y_2$ iff\/ $\cl_{X_2}(f(X_1))=X_2$;

\smallskip

\noindent(e) $g$ is an injection iff $f$ satisfies the following condition:

\smallskip

\noindent{\rm(ZI)} For every $F_1,F_2\in I_1$ such that $F_1\cap F_2=\ems$
there exist $G_1,G_2\in I_2$ with $G_1\cap G_2=\ems$ and $f(F_i)\sbe G_i$, $i=1,2$;

  \smallskip

\noindent(f) $g$ is an open injection iff $I_1\sbe f\inv(I_2)$ and $f$ satisfies condition {\rm (ZO)};

  \smallskip

\noindent(g) $g$ is a closed injection iff $f\inv(I_2)=I_1$;

  \smallskip

\noindent(h) $g$ is a perfect surjection iff $f$ satisfies condition
{\rm (ZP)} and $\cl_{X_2}(f(X_1))=X_2$;

\smallskip

\noindent(i) $g$ is a dense embedding iff\/  $\cl_{X_2}(f(X_1))=X_2$ and $I_1\sbe f\inv(I_2)$.
\end{theorem}

\doc Set $\psi_g=\TE^t(g)$ (see Theorem \ref{genstonec}). Then
$\psi_g:CO(Y_2)\lra CO(Y_1)$, $G\mapsto g\inv(G)$. Set also
$\psi_f: A_2\lra A_1$, $G\mapsto f\inv(G)$. Then, (\ref{0ic}),
(\ref{0i}) and (\ref{psif}) imply that
$\psi_f=r_1^c\circ\psi_g\circ e_2^c$.

\smallskip

\noindent(a) It follows from Lemma \ref{skelnewnew}.

  \smallskip

\noindent(b) {\em First Proof.}~ Using \cite[Theorem 2.8(a)]{Di5}
and (\ref{0il}), we get that the map $g$ is open iff there exists
a map $\psi^f:I_1\lra I_2$ satisfying the following conditions:

\smallskip

\noindent(OZL1) For every $F\in I_1$ and every $G\in I_2$,
$(F\cap f\inv(G)=\ems)\rightarrow (\psi^f(F)\cap G=\ems)$;

\smallskip

\noindent(OZL2) For every $F\in I_1$, $f\inv(\psi^f(F))\spe F$.

Obviously, condition (OZL2) is equivalent to the following one:
for every $F\in I_1$, $f(F)\sbe\psi^f(F)$. We will show that for
every $F\in I_1$, $\psi^f(F)\sbe\cl_{X_2}(f(F))$. Indeed, let
$y\in\psi^f(F)$ and suppose that $y\nin \cl_{X_2}(f(F))$. Since
$I_2$ is a base of $X_2$, there exists a $G\in I_2$ such that
$y\in G$ and $G\cap f(F)=\ems$. Then $F\cap f\inv(G)=\ems$ and
condition (OZL1) implies that $\psi^f(F)\cap G=\ems$. We get that
$y\nin\psi^f(F)$, a contradiction. Thus
$f(F)\sbe\psi^f(F)\sbe\cl_{X_2}(f(F))$. Since $\psi^f(F)$ is a
closed set, we obtain that $\psi^f(F)=\cl_{X_2}(f(F))$. Obviously,
conditions (OZL1) and (OZL2) are satisfied when
$\psi^f(F)=\cl_{X_2}(f(F))$. This implies that $g$ is an open map
iff for every $F\in I_1$, $\cl_{X_2}(f(F))\in I_2$.

\smallskip

\noindent{\em Second Proof.}~ We have, by (\ref{01}), that
$I_i=(f_i)\inv(CK(Y_i))$, for $i=1,2$. Thus, for every $F\in I_i$,
where $i\in\{1,2\}$, we have that $\cl_{Y_i}(f_i(F))\in CK(Y_i)$.

Let $g$ be an open map and $F\in I_1$. Then,
$G=\cl_{Y_1}(f_1(F))\in CK(Y_1)$. Thus $g(G)\in CK(Y_2)$. Since
$G$ is compact, we have that
$g(G)=\cl_{Y_2}(g(f_1(F)))=\cl_{Y_2}(f_2(f(F)))=\cl_{Y_2}(f_2(\cl_{X_2}(f(F))))$.
Therefore, $\cl_{X_2}(f(F))=(f_2)\inv(g(G))$, i.e.
$\cl_{X_2}(f(F))\in I_2$.

Conversely, let $f$ satisfies condition (ZO). Since $CK(Y_1)$ is
an open base of $Y_1$, for showing that $g$ is an open map, it is
enough to prove that for every $G\in CK(Y_1)$,
$g(G)=\cl_{Y_2}(f_2(\cl_{X_2}(f(F))))$ holds, where
$F=(f_1)\inv(G)$ and thus $F\in I_1$. Obviously,
$G=\cl_{Y_1}(f_1(F))$. Using again the fact that $G$ is compact,
we get that $g(G)=g(\cl_{Y_1}(f_1(F)))=\cl_{Y_2}(g(f_1(F)))=
\cl_{Y_2}(f_2(f(F)))=\cl_{Y_2}(f_2(\cl_{X_2}(f(F))))$. So, $g$ is
an open map.

\smallskip

\noindent(c) Since $Y_2$ is a locally compact Hausdorff space and
$CK(Y_2)$ is a base of $Y_2$, we get, using the well-known
\cite[Theorem 3.7.18]{E2}, that $g$ is a perfect map iff
$g\inv(G)\in CK(Y_1)$ for every $G\in CK(Y_2)$. Thus $g$ is a
perfect map iff  $\psi_g(G)\in CK(Y_1)$ for every $G\in CK(Y_2)$.
Now, (\ref{0il}) and (\ref{0ic}) imply that $g$ is a perfect map
$\iff$ $\psi_f(G)\in I_1$ for every $G\in I_2$ $\iff$ $f$
satisfies condition (ZP).

\smallskip

\noindent(d) This is obvious.

\smallskip

\noindent(e) Having in mind (\ref{0il}) and (\ref{0ic}),
our assertion follows from \cite[Theorem 3.5]{Di5}.

\smallskip

\noindent(f) It follows from (b), (\ref{0il}), (\ref{0ic}), and
\cite[Theorem 3.12]{Di5}.

\smallskip

\noindent(g) It follows from (c), (\ref{0il}), (\ref{0ic}), and
\cite[Theorem 3.14]{Di5}.

\smallskip

\noindent(h) It follows from (c) and (d).

\smallskip

\noindent(i) It follows from (d) and \cite[Theorem 3.28 and
Proposition 3.3]{Di5}. We will also give  a {\em second proof}\/
of this fact. Obviously, if $g$ is a dense embedding then $g(Y_1)$
is an open subset of $Y_2$ (because $Y_1$ is locally compact);
thus $g$ is an open mapping and we can apply (f) and (d).
Conversely, if $\cl_{X_2}(f(X_1))=X_2$ and $I_1\sbe f\inv(I_2)$,
then, by (d), $g(Y_1)$ is a dense subset of $Y_2$. We will show
that $f$ satisfies condition (ZO). Let $F_1\in I_1$. Then there
exists $F_2\in I_2$ such that $F_1=f\inv(F_2)$. Then, obviously,
$\cl_{X_2}(f(F_1))\sbe F_2$. Suppose that
$G_2=F_2\stm\cl_{X_2}(f(F_1))\nes$. Since $G_2$ is open, there
exists $x_2\in G_2\cap f(X_1)$. Then there exists $x_1\in X_1$
such that $f(x_1)=x_2\in F_2$. Thus $x_1\in F_1$, a contradiction.
Therefore, $\cl_{X_2}(f(F_1))= F_2$. Thus, $\cl_{X_2}(f(F_1))\in
I_2$. So, condition (ZO) is fulfilled. Hence, by (b), $g$ is an
open map. Now, using (f), we get that $g$ is also an injection.
All this shows that $g$ is a dense embedding.  \sqs

Recall that a continuous map $f:X\lra Y$ is called {\em
quasi-open\/} (\cite{MP}) if for every non-empty open subset $U$
of $X$, $\int(f(U))\nes$ holds. As it is shown in \cite{D1}, if
$X$ is regular and Hausdorff, and $f:X\lra Y$ is a closed map,
then $f$ is quasi-open iff $f$ is skeletal. This fact and Theorem
\ref{zdextcmain} imply the following two corollaries:

\begin{cor}\label{zdextcmaincb}
Let $X_1$, $X_2$ be two zero-dimensional Hausdorff spaces and
  $f:X_1\lra X_2$ be a continuous function. Then:

  \smallskip

\noindent(a) $\b_0f$ is quasi-open iff $f$ is skeletal;

  \smallskip

\noindent(b) $\b_0f$ is an open map iff $f$ satisfies the following condition:

\smallskip

\noindent{\rm(ZOB)} For every $F\in CO(X_1)$,  $\cl_{X_2}(f(F))\in CO(X_2)$ holds;

\smallskip

\noindent(c) $\b_0f$ is a surjection iff\/ $\cl_{X_2}(f(X_1))=X_2$;

 \smallskip

\noindent(d) $\b_0f$ is an injection iff $f\inv(CO(X_2))=CO(X_1)$.
\end{cor}

\begin{cor}\label{zdextcmaincc}
Let $X_1$, $X_2$ be two zero-dimensional Hausdorff spaces,
  $f:X_1\lra X_2$ be a continuous function, $\BB$ be a Boolean algebra admissible for $X_2$, $(cX_2,c)$ be the Hausdorff zero-dimensional
  compactification of $X_2$ corresponding to $\BB$ (see Theorems \ref{dwingerlc} and \ref{dwinger}) and $g:\b_0X_1\lra cX_2$ be the continuous function such
  that $g\circ \b_0=c\circ f$ (its existence is guaranteed by
  Theorem \ref{zdextcb}). Then:

  \smallskip

\noindent(a) $g$ is quasi-open iff $f$ is skeletal;

  \smallskip

\noindent(b) $g$ is an open map iff $f$ satisfies the following condition:

\smallskip

\noindent{\rm(ZOC)} For every $F\in CO(X_1)$,  $\cl_{X_2}(f(F))\in \BB$ holds;

\smallskip

\noindent(c) $g$ is a surjection iff\/ $\cl_{X_2}(f(X_1))=X_2$;

 \smallskip

\noindent(d) $g$ is an injection iff $f\inv(\BB)=CO(X_1)$.
\end{cor}

\end{document}